%% file: SSci.tex
\documentclass[sts]{imsart}

\RequirePackage{amsthm,amsmath,amsfonts,amssymb}
\RequirePackage[numbers]{natbib}

\startlocaldefs

\usepackage{amsfonts, amsmath, amssymb, bbm}

\usepackage{times}
\usepackage{bm}
\usepackage{moreverb,url}
\usepackage{multirow}
\usepackage{float}
\usepackage{tikz}

\usepackage{mathtools}
\mathtoolsset{showonlyrefs}
\graphicspath{{../Figures/}}
\usepackage{color}
\definecolor{darkred}{RGB}{100,0,0}
\definecolor{darkgreen}{RGB}{0,100,0}
\definecolor{darkblue}{RGB}{0,0,150}

\definecolor{purple}{rgb}{0.4,.1,.9}

\usepackage{prodint}

\usepackage{hyperref}
\hypersetup{colorlinks=true, linkcolor=darkred, citecolor=darkgreen, urlcolor=darkblue}

\usepackage[plain,noend]{algorithm2e}

\input{notation}

\endlocaldefs

\begin{document}

\begin{frontmatter}

\title{Asymptotic Theory for Doubly Robust Estimators with Continuous-Time Nuisance Parameters}
\runtitle{Asymptotic Theory for Doubly Robust Estimators with Continuous-Time Nuisance Parameters}

\begin{aug}
\author[A]{\fnms{Andrew}~\snm{Ying}\ead[label=e1]{aying9339@gmail.com}}
\address[A]{Andrew Ying\printead[presep={\ }]{e1}.}
\end{aug}

\begin{abstract}
Doubly robust estimators have gained widespread popularity in various fields due to their ability to provide unbiased estimates under model misspecification. However, the asymptotic theory for doubly robust estimators with continuous-time nuisance parameters remains largely unexplored. In this short communication, we address this gap by developing a general asymptotic theory for a class of doubly robust estimating equations involving stochastic processes and Riemann-Stieltjes integrals. We introduce generic assumptions on the nuisance parameter estimators that ensure the consistency and asymptotic normality of the resulting doubly robust estimator. Our results cover both the model doubly robust estimator, which relies on parametric or semiparametric models, and the rate doubly robust estimator, which allows for flexible machine learning methods. We discuss the implications of our findings and highlight the key differences between the continuous-time setting and the classical theory for doubly robust estimators. Our work provides a solid theoretical foundation for the use of doubly robust estimators in complex settings with continuous-time nuisance parameters, paving the way for future research and applications.
\end{abstract}



\end{frontmatter}



\section{Introduction}
Doubly robust estimators have gained significant attention in recent years due to their attractive property of providing unbiased estimates even when one of the two nuisance models is misspecified. This robustness against model misspecification has led to their widespread use in various fields, such as missing data analysis \citep{bang2005doubly,  tsiatis2006semiparametric} and causal inference \citep{robins1995analysis, tchetgen2010doubly,  van2003unified, ghassami2022minimax}. The classical theory for doubly robust estimators has been well-established, particularly for settings without complex time-dependent structures \citep{ying2022proximal, wang2022doubly}.

However, in many applications, the nuisance parameters involved in the estimating equations may take the form of stochastic processes, requiring the use of Riemann-Stieltjes integrals. Examples of such settings include informative censoring in survival analysis \citep{rotnitzky2005inverse, ying2022proximal}, covariate-induced dependent truncation \cite{wang2022doubly}, and irregular longitudinal data \citep{rytgaard2022continuous}. Despite the growing need for doubly robust estimators in these complex scenarios, the asymptotic theory for such estimators with continuous-time nuisance parameters remains largely unexplored.

In this short communication, we aim to fill this gap by developing a general asymptotic theory for a class of doubly robust estimating equations involving stochastic processes and Riemann-Stieltjes integrals. We introduce generic assumptions on the nuisance parameter estimators that ensure the consistency and asymptotic normality of the resulting doubly robust estimator. Our results cover both the model doubly robust estimator, which relies on parametric or semiparametric models, and the rate doubly robust estimator, which allows for flexible machine learning methods.

A key challenge in establishing the asymptotic properties of doubly robust estimators with continuous-time nuisance parameters lies in the need for stronger assumptions on the total variation of the estimated processes. We discuss these assumptions and their implications, highlighting the differences between the continuous-time setting and the classical theory for doubly robust estimators.

The remainder of this paper is organized as follows. Section \ref{sec:doubly} introduces the general class of continuous-time doubly robust estimating equations and provides examples of their applications. Section \ref{sec:ral} presents the asymptotic theory for model double robustness, while Section \ref{sec:cross} focuses on the asymptotic theory for rate double robustness. We conclude the paper in Section \ref{sec:dis} with a summary of our findings and a discussion of potential future directions.

\section{Continuous-Time Doubly Robust Estimating Equations}\label{sec:doubly}

\subsection{General Class of Estimating Equations}

Consider two stochastic processes $H(t; \theta)$ and $Q(t; \theta)$ of data $\cO$, where $\theta$ is the parameter of interest. We consider the following general class of estimating functions
\begin{align}
    \Xi\{H(\theta), Q\} &= H(\infty; \theta)Q(\infty) - \int^\infty H(t; \theta)dQ(t),\label{eq:continuousform}
\end{align}
where $\theta$ is a finite-dimensional parameter, $H$ and $Q$ are stochastic processes, and the form of $H(\infty; \theta)$ is known to the researcher as initial condition. Here the integral is understood as pathwise Riemann-Stieltjes integration. Denote $H$, $Q$, and $\theta$ as the truth, $\Xi\{H(\theta), Q\}$ satisfies population double robustness
\begin{equation}
    \E[\Xi\{H'(\theta), Q\}] = \E[\Xi\{H(\theta), Q'\}] = 0,
\end{equation}
for any $H'$ and $Q'$.
We impose here a standard population-level assumption, which is needed in the classical case even without nuisance parameters.
\begin{assump}[Donsker]\label{assump:donsker}
\begin{equation}
    \{\Xi\{H(\theta'), Q\} - \Xi\{H(\theta), Q\}: \|\theta' - \theta\| < \delta\}
\end{equation}
is $\P$-Donsker for some $\delta > 0$ and that 
\begin{equation}
    \P[\Xi\{H(\theta'), Q\} - \Xi\{H(\theta), Q\}^2] \to 0,
\end{equation}
as $\theta' \to \theta$.
\end{assump}
This assumption \citep[Lemma 3.3.5]{van1997weak} is typically satisfied for identifying low-dimensional parameters.

\subsection{Examples and Importance}
Here we give some examples to illustrate the generality of our estimating functions.
\begin{exm}[Informative Right Censoring]\label{exm:rightcensoring}

By considering right censoring as monotone coarsening, \citet{van2003unified, rotnitzky2005inverse, tsiatis2006semiparametric} developed augmented inverse probability weighting estimators for survival analysis under the coarsening at random assumption. Define $T$ as event time, $C$ as censoring time, $\Delta = \mathbbm{1}(T \leq C)$ as event indicator. We also define a covariate process $L(t)$ and write $\bar L(t) = \{L(s): s \leq t\}$. We define
\begin{align*}
    &\Xi\{H(\theta), Q\} \\
    &= 
    \frac{\Delta D\{T, \overline L(T); \theta\}}{\Prodi^T \left[1 - \lambda_C\left\{t|\overline L(t)\right\}dt\right]} \\
    &~~+ \int^\infty \frac{\E\left[D\{T, \overline L(T); \theta\}|T \geq t, \overline L(t)\right]}{\Prodi^t \left[1 - \lambda_C\left\{u|\overline L(u)\right\}du\right]}dM_C\{t, \overline L(t)\},
\end{align*}
where
\begin{equation}
    H(t; \theta) = \mathbbm{1}(T \geq t)\E\left[D\{T, \overline L(T); \theta\}|T \geq t, \overline L(t)\right],
\end{equation}
\begin{equation}
    Q(t) = \frac{\mathbbm{1}(C > t)}{\Prodi_{0 \leq u \leq t} \left[1 - \lambda_C\left\{u|\overline L(u)\right\}du\right]}.
\end{equation}
\end{exm}

\begin{exm}[Informative Right Censoring with Proxies]
Unlike \citet{van2003unified, rotnitzky2005inverse, tsiatis2006semiparametric}, though considering right censoring as monotone coarsening, \citet{ying2023proximal} developed doubly estimators for survival analysis without coarsening at random assumption. Following the notation in example \ref{exm:rightcensoring}, he further broke down $L(t)$ into $(X(t), Z(t), W(t)$. Define
\begin{align*}
    &\Xi\{H(\theta), Q\} \\
    &= \Delta \text{BP}_2(\widetilde T) D\{\widetilde T, \overline X(\widetilde T); \theta\}\\
    &~~- \int^{\widetilde T}\text{BP}_1(t; D, \theta)\{d\text{BP}_2(t) - dN_C(t)\text{BP}_2(t)\},
\end{align*}
where
\begin{equation}
    H(t; \theta) = \mathbbm{1}(T \geq t)\text{BP}_1(t),
\end{equation}
\begin{equation}
    Q(t) = \mathbbm{1}(C > t)\text{BP}_2(t),
\end{equation}
and $\text{BP}_1(t; D, \theta) = \text{BP}_1\{t, \overline W(t), \overline X(t); D, \theta\}$ satisfies 
\begin{align*}
    &\E\{d\text{BP}_1(t; D, \theta) - \text{BP}_1(t; D, \theta)dN_T(t)|\widetilde T \geq t, \overline Z(t), \overline X(t)\} \\
    &= 0,
\end{align*}
with an initial condition
\begin{equation}\label{eq:eventbridgeUinit}
    \text{BP}_1(\infty; D, \theta) = D\{T, \overline X(T); \theta\},
\end{equation}
$\text{BP}_2(t) = \text{BP}_2\{t, \overline Z(t), \overline X(t)\}$ satisfies
\begin{equation}\label{eq:censorbridgeU}
    \E\{d\text{BP}_2(t) - \text{BP}_2(t)dN_C(t)|\widetilde T \geq t, \overline w(t), \overline X(t)\} = 0,
\end{equation}
with an initial condition
\begin{equation}\label{eq:censorbridgeUinit}
    \text{BP}_2\{0, \overline z(0), \overline x(0)\} = 1.
\end{equation}

\end{exm}

\begin{exm}[Conditional Quasi-independent Left Truncation]

\citet{wang2022doubly} developed a doubly robust estimator in the presence of covariate-induced dependent left truncation.
\begin{align*}
    &\Xi\{H(\theta), Q\}\\
    &=\frac{D(T, Z; \theta)}{G(T|Z)} \\
    &~~- \int^\infty \frac{\E_F\{\mathbbm{1}(T \leq t)D(T, Z; \theta)\}}{1 - F(t|Z)} \frac{d\bar M_Q(t;G)}{G(t|Z)},
\end{align*}
where
\begin{equation}
    H(t; \theta) = \mathbbm{1}(T > t)\frac{\E_F\{\mathbbm{1}(T \leq t)D(T, Z; \theta)\}}{ 1-F(t|Z) },
\end{equation}
\begin{equation}
    Q(t) = \frac{\mathbbm{1}(Q \leq t)}{G(t|Z)}.
\end{equation}
\end{exm}

\begin{exm}[Continuous-time Causal Inference]
Adhering to the notation introduced by \citet{rytgaard2022continuous}, we define $A(t)$ and $L(t)$ as the marked point processes representing treatment and confounders, where $t \in [0, \tau]$. Consider an end-of-study outcome $Y$ measure at $\tau$ and define $Y_{\bar a}$ as the potential outcome if a patient were treated following $\bar a$. Consider a parameter of interest as $\theta = \E(Y_{\bar 1})$ the average of potential outcome under always-treated. Define the observed data likelihood up to time $t$
\begin{align*}
    &d\P_t \\
    &= \Prodi_{t}d\P(A(s)|\bar A(s-), \bar L(s-))d\P(L(s)|\bar A(s-), \bar L(s-)),
\end{align*}
and the target data likelihood up to time $t$ as
\begin{align*}
    d\P_{\bar 1, t} = \Prodi_{t}d\mathbbm{1}(A(s) = 1)d\P(L(s)|\bar A(s-), \bar L(s-)).
\end{align*}

By \citet{rytgaard2022continuous}, an estimating function is
\begin{align*}
    &\Xi\{H(\theta), Q\} = (Y - \theta)\\
    &~~- \int\int y d\P_{\bar 1, \tau}\{y, \bar a, \bar l|\bar A(t), \bar L(t)\}d\left(\frac{d\P_{\bar 1, t}}{d\P_t}\right),
\end{align*}
where
\begin{equation}
    H(t; \theta) = \int y d\P_{\bar 1, \tau}\{y, \bar a, \bar l|\bar A(t), \bar L(t)\},
\end{equation}
\begin{equation}
    Q(t) = \frac{d\P_{\bar 1, t}}{d\P_t}.
\end{equation}
\end{exm}

\subsection{Why is there a gap intuitively?}\label{sec:gap}

The double robustness without time structure or with discrete time structure are both special cases of \eqref{eq:continuousform}. For example, without time, \eqref{eq:continuousform} becomes
\begin{align*}
    \Xi\{H(\theta), Q\} &= H(1; \theta)Q(0) - H(0; \theta)Q(0) + H(0; \theta).
\end{align*}
The bias when plugging in $H'$ and $Q'$ is
\begin{equation}
    \E\{\Xi(H'(\theta), Q')\} = \E[\{H'(\theta) - H(\theta)\}(Q' - Q)],
\end{equation}
which is a saddle surface (or a hyperbolic paraboloid) of $H$ and $Q$, or known as mixed bias property. Clearly from this, one can tell that the mixed bias property is a sufficient but not necessary condition of population double robustness. However, the rate double robustness is mainly if ever, build on the mixed bias property by applying the Cauchy-Schwartz equality as
\begin{equation}
    |\E\{(H' - H)(Q - Q)\}| = \|(H' - H)\|_2\|Q' - Q\|_2.
\end{equation}
Therefore the bias is controlled by the product error of $H'$ and $Q'$.

When both $H$ and $Q$ involve time and $\Xi(H, Q)$ constitutes a Riemann–Stieltjes integral, applying the Cauchy-Schwartz inequality yields bound with a total variation norm on either $H$ or $Q$. However, a total variation bound on either $H$ or $Q$ can be unreasonable because, simply considering the simplest case when $H$ or $Q$ is the empirical cumulative distribution function, the total variation norm is always of constant rate and cannot vanish over sample size $n$. On the other hand, step functions like empirical cumulative distribution function and nonparametric maximum likelihood estimators are largely employed in continuous-time setting. Therefore, extra caution is warranted.

\section{Model Double Robustness}\label{sec:ral}
In this section, we present the asymptotic theory for model doubly robust estimators in the continuous-time setting. We first introduce the necessary assumptions and then state the main result on the consistency and asymptotic normality of the estimator.
\subsection{Assumptions}
Suppose one is able to model $H$ and $Q$ parametrically or semiparametrically. For instance, in Example \ref{exm:rightcensoring}, $H$ or $Q$ can be estimated once estimating the distribution of $T$ given $\bar L$ or $C$ given $\bar L$. This can be done using Cox proportional hazard model or Aalen additive model, which also implies model on $H$ and $Q$. Both are reasonable and practical semiparametric model that are well studied in the literature. In these cases, the estimators $\hat H$ and $\hat Q$ are asymptotically linear and uniformly root-$n$ consistent.

Denote the i.i.d. observed data $\{\cO_i = \cO_i(t): t \geq 0\}_{i = 1}^n$. We write $H_i(\theta, t)$ and $Q_i(t)$ when plugging in $\cO_i$, and 
\begin{align}
    \Xi_i\{H(\theta), Q\} &= H_i(\infty; \theta)Q_i(\infty) - \int^\infty H_i(t; \theta)dQ_i(t).
\end{align}

We define $\hat \theta_{\text{MDR}}$ as the solution to the estimating equation
\begin{equation}
    \frac{1}{n}\sum_{i = 1}^n\Xi_i\{\hat H(\theta), \hat Q\} = o_p(1/\sqrt{n}).
\end{equation}
$o(1/\sqrt{n})$ is imposed here to accommodate the case when actual zero cannot be achieved, like median. Later we give sufficient and generic conditions on $\hat H$ and $\hat Q$ under which $\hat \theta_{\text{MDR}}$ is consistent and asymptotically normal. We additionally define $\|\cdot\|_2$ as the $L^2$-norm. Additionally, for any stochastic $f(\cO, t)$ we define the $L^2$ supremum norm
\begin{equation}
    \|f\|_{\sup, 2} = \left\|\sup_t|f(\cO, t)|\right\|_2,
\end{equation}
and $L^2$ total variation norm
\begin{align*}
    &\|f\|_{\text{TV}, 2} \\
    &= \left\||f(\cO, 0)| +\sup_{\cP}\sum_{m = 0}^M|f(\cO, t_{m + 1}) - f(\cO, t_m)|\right\|_2.
\end{align*}
where the supremum runs over the set of all partitions $\cP =\left\{P = \{t,\cdots, t_{M}\}\mid P{\text{ is a partition of }}[0,\infty)\right\}$ of the given interval.

\begin{assump}[Uniform convergence]\label{assump:uniformcons1}
~~~
\begin{enumerate}
    \item We have
\begin{equation}
    \left\|\hat H(\theta) - H_*(\theta)\right\|_{\sup, 2} = o(1),
\end{equation}
for some $H_*(t; \theta)$ with $\left\|H_*(\theta)\right\|_{\text{TV}, 2} < \infty$.
    \item We have
\begin{equation}
    \left\|\hat Q - Q_*\right\|_{\sup, 2} = o(1),
\end{equation}
for some $Q_*(t)$ with $\left\|Q_*\right\|_{\text{TV}, 2} < \infty$.
\end{enumerate}
\end{assump}
The uniform convergence assumption requires that our estimators to converge uniformly to some limiting processes in probability, but not necessarily to the truths. It is usually satisfied for estimators based on Gilvenko-Cantelli class \citep{van1997weak}. The following asymptotic linearity assumption is needed to ensure asymptotic normality.
\begin{assump}[Asymptotic linearity]\label{assump:if}
For some limiting processes $H_*(\cO; t)$ and $Q_*(\cO; t)$, and some mean zero i.i.d. stochastic processes $\xi_{h, i}(\cO, t)$ and $\xi_{q, i}(\cO, t)$, we have either one of the following holds:
{\small
\begin{itemize}
    \item 
\begin{equation}
    \left\|\hat H(\cO; \theta) - H_*(\cO; \theta) - \frac{1}{n}\sum_{i = 1}^n \xi_{h, i}(\cO)\right\|_{\sup, 2} = o(n^{-\frac{1}{2}}),
\end{equation}
and
\begin{equation}
    \left\|\hat Q(\cO) - Q_*(\cO) - \frac{1}{n}\sum_{i = 1}^n \xi_{q, i}(\cO)\right\|_{\text{TV}, 2} = o(n^{-\frac{1}{2}});
\end{equation}
\item 
\begin{equation}
    \left\|\hat H(\cO; \theta) - H_*(\cO; \theta) - \frac{1}{n}\sum_{i = 1}^n \xi_{h, i}(\cO)\right\|_{\text{TV}, 2} = o(n^{-\frac{1}{2}}),
\end{equation}
and
\begin{equation}
    \left\|\hat Q(\cO) - Q_*(\cO) - \frac{1}{n}\sum_{i = 1}^n \xi_{q, i}(\cO)\right\|_{\sup, 2} = o(n^{-\frac{1}{2}}).
\end{equation}
\end{itemize}
}

\end{assump}
The asymptotic linearity assumption requires both estimators can be asymptotically linear expanded to $L_2$ supremum norm and at least one of the expansion can be strengthened into total variation norm. It is usually satisfied with estimators based on Donsker class \citep{van1997weak}. \citet{wang2022doubly} have shown that estimated conditional cumulative hazard functions output from Cox regression and Aalen regression satisfy Assumptions \ref{assump:uniformcons1}, \ref{assump:if}. This type of asymptotic linearity in the total variation sense was also considered in \citet{hou2021treatment, ying2021proximal}.

\subsection{Main Result}
Under the assumptions stated in Section 3.1, we establish the consistency and asymptotic normality of the model doubly robust estimator $\hat \theta_{\text{MDR}}$, which is defined as the solution to the estimating equation:

\begin{thm}\label{thm:mdr}
Under Assumptions \ref{assump:donsker} and \ref{assump:uniformcons1}, and either $H_* = H$ or $Q_* = Q$, then
\begin{equation}
    \hat \theta_{\text{MDR}} \to \theta,
\end{equation}
in probability and further if Assumption \ref{assump:if} holds, then
\begin{align}
    &\left[\frac{\partial}{\partial \theta}\E\{\Xi(H_*(\theta), Q_*\}\right]\cdot\sqrt{n}(\hat \theta_{\text{MDR}} - \theta) \\
    &\to \cN(0, \E(\Xi\{H_*(\theta), Q_*\}^2\}),
\end{align}
in distribution.
\end{thm}

Theorem \ref{thm:mdr} establishes the consistency and asymptotic normality of the model doubly robust estimator when one of the nuisance processes is consistently estimated using parametric or semiparametric models. 
\subsection{Discussion and Implications}

The uniform convergence assumption (Assumption \ref{assump:uniformcons1}) ensures that the estimators $\hat H(t; \theta)$ and $\hat Q(t)$ converge uniformly to their respective limits, while the asymptotic linearity assumption (Assumption \ref{thm:mdr}) guarantees that at least one of the estimators satisfies the necessary asymptotic expansion condition. These assumptions are crucial for establishing the consistency and asymptotic normality of the model doubly robust estimator.

Compared to the classical Z-estimation theory, the continuous-time setting requires stronger assumptions on the total variation of the estimated processes. This is due to the presence of Riemann-Stieltjes integrals in the estimating equations, which necessitate the control of the total variation of the nuisance process estimators.

The implications of Theorem \ref{thm:mdr} are significant for the use of model doubly robust estimators in practice. It provides a rigorous theoretical foundation for the consistency and asymptotic normality of these estimators in the continuous-time setting, enabling valid inference and confidence interval construction.

However, it is important to note that the model double robustness property relies on the correct specification of at least one of the nuisance processes. In practice, it may be challenging to ensure the correct specification of parametric or semiparametric models for the nuisance processes, especially in the presence of complex time-dependent structures. This motivates the development of rate doubly robust estimators below, which allow for more flexible estimation of the nuisance processes using machine learning methods.

\section{Rate Double Robustness}\label{sec:cross}
In this section, we present the asymptotic theory for rate doubly robust estimators in the continuous-time setting. We first introduce the estimator construction and the cross-fitting procedure, a typical procedure applied to this kind of estimators \citep{chernozhukov2018double}, followed by the necessary assumptions and the main result on the consistency and asymptotic normality of the estimator.

\subsection{Estimator Construction and Cross-Fitting}
The estimator is defined as follows. We randomly split $\{1, \cdots, n\}$ into $L$ folds of equal size with the index sets $\cI_1,\cdots,\cI_L$. For $1 \leq l \leq L$, estimate the nuisance processes $H(t)$ and $Q(t)$ using the out-of-l-fold data indexed by $\cI_{-l} = \{1, \cdots, n\}/\cI_l$, and denote the estimates by $\hat H^{(-l)}(\theta)$ and $\hat Q^{(-l)}$.

We define $\hat \theta_{\text{RDR}}$ as the solution to the estimating equation
\begin{equation}\label{eq:dree}
    \frac{1}{n}\sum_{l = 1}^L\sum_{i = 1}^n\Xi_i\{\hat H^{(-l)}(\theta), \hat Q^{(-l)}\} = o_p(1/\sqrt{n}).
\end{equation}
Researchers aim to impose assumptions under which they can obtain root-$n$ asymptotically normal estimators. 

\subsection{Assumptions}
To establish the asymptotic properties of the rate doubly robust estimator, we require the following assumptions:
\begin{assump}[Uniform consistency]\label{assump:uniformcons2}
The uniform $L_2$ norms converge to zero, that is,
\begin{equation}
    \sup_{\theta}\E_{\cO_{\cI_1}}\left[
    \E_{\cO_{\cI_{-1}}}\left\{\sup_t|\hat H(\theta) - H(\theta)|^2\right\}\right]^{1/2} = o(1).
\end{equation}
\begin{equation}
    \E_{\cO_{\cI_1}}\left[
    \E_{\cO_{\cI_{-1}}}\left\{\sup_t|\hat Q - Q|^2\right\}\right]^{1/2} = o(1).
\end{equation}
\end{assump}
The uniform consistency requires that our estimators converge uniformly to a limit process in probability. The supremum is over both time and covariates process. This is usually satisfied by most machine learners. 

\begin{assump}[Rate condition]\label{assump:rate}
The uniform $L_2$ norm and the $L_2$ total variation norm both converge to zero at some speed, for any $\E\{Q_n(t)\} = O(1)$,
\begin{align}
    &\E_{\cO_{\cI_1}}\left(
    \E_{\cO_{\cI_{-1}}}\left[\int\{\hat H(t; \theta) - H(t; \theta)\}d\{\hat Q(t) - Q(t)\}\right]\right) \\
    &= o(1/\sqrt{n}).\label{eq:cross}
\end{align}
\end{assump}
These assumptions ensure that the cross integrated error term \eqref{eq:cross} converges faster than root-$n$. This is a unique assumption in the continuous-time literature compared to non-continuous-time cases. In non-continuous-time cases, Assumption \ref{assump:rate} is usually imposed as the product error rate of $\hat H$ and $\hat Q$ being faster than root-$n$. However, in continuous-time cases, we have shown that this leads to an unrealistic total variation norm. Sufficient conditions of Assumption \ref{assump:rate} is interesting. For example, \citet{rytgaard2022continuous} considered smooth estimators $\hat H$ and $\hat Q$ to satisfy this assumption.
\subsection{Main Result}
Under the assumptions stated in Sections 3.1 and 4.2, we establish the consistency and asymptotic normality of the rate doubly robust estimator $\theta_{\text{RDR}}$.
\begin{thm}\label{thm:rdr}
Under Assumptions \ref{assump:donsker}, \ref{assump:uniformcons2}, 
\begin{equation}
    \hat \theta_{\text{RDR}} \to \theta,
\end{equation}
in probability and furthermore if Assumption \ref{assump:rate} holds, then
\begin{align}
    &\left[\frac{\partial}{\partial \theta}\E\{\Xi(H(\theta), Q\}\right]\cdot\sqrt{n}(\hat \theta_{\text{RDR}} - \theta) \\
    &\to \cN(0, \E[\Xi\{H(\theta), Q\}^2]),
\end{align}
in distribution.
\end{thm}

Theorem \ref{thm:rdr} establishes the consistency and asymptotic normality of the rate doubly robust estimator when the product of the convergence rates of the nuisance process estimators is faster than the parametric rate $n^{-1/2}$. This condition allows for the use of flexible machine learning methods to estimate the nuisance processes, as long as their combined convergence rate is sufficiently fast.

The proof of Theorem \ref{thm:rdr} relies on the decomposition of the estimating equation and the control of the resulting terms using the assumptions on the nuisance process estimators. The cross-fitting procedure plays a crucial role in ensuring that the bias terms are asymptotically negligible.

\section{Discussion}\label{sec:dis}
In this paper, we have developed a general asymptotic theory for doubly robust estimators with continuous-time nuisance parameters. We considered a broad class of estimating equations involving stochastic processes and Riemann-Stieltjes integrals, which encompass a wide range of applications in various fields, such as survival analysis and causal inference. We established the consistency and asymptotic normality of the model doubly robust estimator under suitable assumptions on the uniform convergence and asymptotic linearity of the nuisance process estimators. These assumptions are specific to the continuous-time setting and differ from those in the classical Z-estimation theory. We also provided a rigorous theoretical foundation for the use of rate doubly robust estimators, which allow for flexible machine learning methods to estimate the nuisance processes, under the condition that their combined convergence rate is faster than the parametric rate.

One of the key insights of our work is the identification of the gap between non-continuous-time double robustness and continuous-time double robustness. We highlighted the challenges in establishing continuous-time double robustness with additional assumptions on the total variation of the estimated processes. This gap underscores the importance of carefully considering the assumptions and convergence rates when using doubly robust estimators in practice.

Our theoretical results have important implications for the development and application of doubly robust estimators in settings with complex time-dependent structures. They provide guidance on the construction of suitable estimators, the choice of appropriate nuisance process models, and the assessment of the estimators' asymptotic properties. The general framework we introduced can serve as a foundation for further research on doubly robust estimation in various domains. In fact, the asymptotic theory proofs in \citet{ying2022proximal, wang2022doubly, luo2022doubly} were inspired by an early version of this work.

In conclusion, our work advances the understanding of doubly robust estimation in the presence of continuous-time nuisance parameters and provides a rigorous theoretical foundation for the use of these estimators in a wide range of applications. We hope that our findings will stimulate further research on the design and analysis of doubly robust estimators in complex settings and contribute to the development of robust and efficient statistical methods for the analysis of time-dependent data.

\begin{acks}[Acknowledgments]
The author would like to thank Yuyao Wang, Ronghui Xu for countless discussions. The author would like to acknowledge the helpful discussions and insights gained through interactions with the Claude AI assistant, which contributed to refining the presentation of this work.
\end{acks}

\appendix
\section{Proofs}

\subsection{Proof of Theorem \ref{thm:mdr}}
We have the following decomposition
\begin{align}
    &\frac{1}{n}\sum_{i = 1}^{n}\Xi_i\{\hat H(t; \theta), \hat Q(t)\}\\
    &= \cT_1 + \cT_2 + \cT_3 + \cT_4 + \cT_5 + \cT_6 + \cT_7,
\end{align}
where $H_*$ and $Q_*$ are some processes (later chosen as the limit process), 
\begin{align*}
    \cT_1 &= \frac{1}{n}\sum_{i = 1}^{n}\Xi_i\{\hat H(\theta), \hat Q\} + \frac{1}{n}\sum_{i = 1}^{n}\Xi_i\{\hat H(\theta), Q_*\}\\
    &~~~- \frac{1}{n}\sum_{i = 1}^{n}\Xi_i\{H_*(\theta), \hat Q\} + \frac{1}{n}\sum_{i = 1}^{n}\Xi_i\{H_*(\theta), Q_*\},\\
    \cT_2 &= \frac{1}{n}\sum_{i = 1}^{n}[\Xi_i\{\hat H(\theta), Q_*\} - \Xi_i\{H_*(\theta), Q_*\}],\\
    \cT_3 &= \frac{1}{n}\sum_{i = 1}^{n}[\Xi_i\{H_*(\theta), \hat Q\} - \Xi_i\{H_*(\theta), Q_*\}],\\
    \cT_4 &= \frac{1}{n}\sum_{i = 1}^{n}(\Xi_i\{H_*(\theta), Q_*\} - \E[\Xi\{H_*(\theta), Q_*\}]),\\
    &~~~- \frac{1}{n}\sum_{i = 1}^{n}(\Xi_i\{H_*(\theta), Q_*\} - \E[\Xi\{H_*(\theta), Q_*\}]),\\
    \cT_5 &= \frac{1}{n}\sum_{i = 1}^{n}(\Xi_i\{H_*(\theta), Q_*\} - \E[\Xi\{H_*(\theta), Q_*\}]),\\
    \cT_6 &= \E[\Xi\{H_*(\theta), Q_*\}] - \E[\Xi\{H(\theta), Q\}], \\
    \cT_7 &= \E[\Xi\{H_*(\theta), Q_*\}] - \E[\Xi\{H(\theta), Q\}].
\end{align*}

We directly prove for asymptotic normality since consistency proof is similar and simpler. By Assumption \ref{assump:if} and Chebyshev's inequality 
\begin{align}
    &\P(\sqrt{n}|\cT_1| > \delta) \\
    &\leq \frac{\delta^2}{n}\Var\left\{\sum_{i = 1}^n \frac{1}{n^2}\int^\infty \sum_{j = 1}^n\xi_j(t; \cO_i)d\sum_{k = 1}^n\xi_k(t; \cO_i)\right\} \\
    &= \frac{\delta^2}{n^5}\sum_{i, i', j, j', k, k' = 1}^n \E\Bigg\{\int^\infty \xi_j(t; \cO_i)d\xi_k(t; \cO_i)\\
    &~~~\int^\infty \xi_{j'}(t; \cO_{i'})d\xi_{k'}(t; \cO_{i'})\Bigg\} \\
    &\leq \frac{\delta^2}{n^5}\sum_{i, i, j, k = 1}^n \E\Bigg\{\int^\infty \xi_j(t; \cO_i)d\xi_k(t;  \cO_i)\\
    &~~~\int^\infty \xi_{j}(t; \cO_{i'})d\xi_{k}(t; \cO_{i'})\Bigg\} \\
    &~~~+ \frac{\delta^2}{n^5}\sum_{i, i', j, k = 1}^n \E\Bigg\{\int^\infty \xi_j(t; \cO_i)d\xi_j(t; \cO_i)\\
    &~~~\int^\infty \xi_k(t; \cO_{i'})d\xi_k(t; \cO_{i'})\Bigg\} \\
    &=O\left(\frac{1}{n}\right) = o(1).
\end{align}
The same arguments with Assumption \ref{assump:if} and Chebyshev's inequality can also show that $\sqrt{n}\cT_2$ and $\sqrt{n}\cT_3$ are $o_p(1)$. $\cT_4$ by Assumption \ref{assump:donsker} and \citet[Lemma 3.3.5]{van1997weak} is $o_p(1 + \sqrt{n}\|\theta - \theta\|)$. $\cT_5$ by weak law of large numbers converges to zero in probability. $\cT_6$ contains our estimator. 

\subsection{Proof of Theorem \ref{thm:rdr}}
We have the following decomposition
\begin{align}
    &\frac{1}{n}\sum_{l = 1}^L\sum_{i \in \cI_l}\Xi_i\{\hat H^{(-l)}(\theta), \hat Q^{(-l)}\}\\
    &= \cT_1 + \cT_2 + \cT_3 + \cT_4 + \cT_5 + \cT_6 + \cT_7,
\end{align}
where 
\begin{align*}
    \cT_1 &= \frac{1}{n}\sum_{l = 1}^L\sum_{i \in \cI_l}\Xi_i\{\hat H^{(-l)}(\theta), \hat Q^{(-l)}\} \\
    &~~~- \frac{1}{n}\sum_{l = 1}^L\sum_{i \in \cI_l}\Xi_i\{\hat H^{(-l)}(\theta), Q\}\\
    &~~~- \frac{1}{n}\sum_{l = 1}^L\sum_{i \in \cI_l}\Xi_i\{H(\theta), \hat Q^{(-l)}\} \\
    &~~~+ \frac{1}{n}\sum_{l = 1}^L\sum_{i \in \cI_l}\Xi_i\{H(\theta), Q\},\\
    \cT_2 &= \frac{1}{n}\sum_{l = 1}^L\sum_{i \in \cI_l}[\Xi_i\{\hat H^{(-l)}(\theta), Q\} - \Xi_i\{H(\theta), Q\}],\\
    \cT_3 &= \frac{1}{n}\sum_{l = 1}^L\sum_{i \in \cI_l}[\Xi_i\{H(\theta), \hat Q^{(-l)}\} - \Xi_i\{H(\theta), Q\}],\\
    \cT_4 &= \frac{1}{n}\sum_{l = 1}^L\sum_{i \in \cI_l}(\Xi_i\{H(\theta), Q\} - \E[\Xi\{H(\theta), Q\}]),\\
    &~~~- \frac{1}{n}\sum_{l = 1}^L\sum_{i \in \cI_l}(\Xi_i\{H(\theta), Q\} - \E[\Xi\{H(\theta), Q\}]),\\
    \cT_5 &= \frac{1}{n}\sum_{l = 1}^L\sum_{i \in \cI_l}(\Xi_i\{H(\theta), Q\} - \E[\Xi\{H(\theta), Q\}]),\\
    \cT_6 &= \E[\Xi\{H(\theta), Q\}] - \E[\Xi\{H(\theta), Q\}], \\
    \cT_7 &= \E[\Xi\{H(\theta), Q\}] - \E[\Xi\{H(\theta), Q\}].
\end{align*}
We directly prove for asymptotic normality since consistency proof is similar and simpler. Firstly, by linearity,
\begin{align*}
    \cT_{1} = \sum_{l = 1}^L\cT_{1, L},
\end{align*}
where
{\footnotesize
\begin{align*}
    \cT_{1, L} = \frac{1}{n}\sum_{i \in \cI_l} \int^\infty \{\hat H_i^{(-l)}(t; \theta) - H_{i}(t; \theta)\}d\{\hat Q_i^{(-l)}(t) -  Q_i(t)\},
\end{align*}
}
whereby Markov's inequality, symmetry across $i$, and Holder's inequality, we have that for any fixed $\delta > 0$,
{\tiny
\begin{align*}
    &\P(\sqrt{n}|\cT_{1, L}| \geq \delta)\\
    &\leq \frac{1}{\delta\sqrt{n}}\E_{\cO_{\cI_{-l}}}\left(\E_{\cO_{\cI_{l}}}\left[\left|\sum_{i \in \cI_l} \int \{\hat H_i(t; \theta) - H_{i}(t; \theta)\}d\{\hat Q_i(t) - Q_{i}(t)\}\right|\right]\right) \\
    &=\frac{\sqrt{n}}{\delta L}\E_{\cO_{\cI_{-l}}}\left(\E_{\cO_{\cI_{l}}}\left[\left|\int \{\hat H(t; \theta) - H(t; \theta)\}d\{\hat Q(t) - Q(t)\}\right|\right]\right) \\
    &= o(1),
\end{align*}
}
by Assumption \ref{assump:rate}. With a symmetric argument, we conclude, by our assumption on the product rate, that
\begin{equation}
    \P(|\cT_{1, L}| \geq \delta) = o(1).
\end{equation}
Since $L$ is constant relatively to $n$, so is $\cT_1$. For $\cT_2$ we write more explicitly as
\begin{equation}
    \cT_2 = \sum_{l = 1}^L\cT_{2, l} 
\end{equation}
where
\begin{align}
    \cT_{2, l} = \frac{1}{n}\sum_{i \in \cI_l} \int^\infty \{\hat H_i^{(-l)}(t; \theta) - H_{i}(t; \theta)\}dQ_{i}(t).
\end{align}
By Chebyshev's inequality, the law of total variance, independence between individuals, and Fatou's lemma we have
{\tiny
\begin{align*}
    &\P(\sqrt{n}|\cT_{2, l}| \geq \delta) \\
    &\leq \frac{1}{n\delta^2} \Var(\cT_{2, l})\\
    &=\frac{1}{n\delta^2} \Var\left[\sum_{i \in \cI_l} \int^\infty \{\hat H_i^{(-l)}(t; \theta) - H_{i}(t; \theta)\}dQ_{i}(t)\right]\\
    &=\frac{1}{n\delta^2}\E\left(\Var\left[\sum_{i \in \cI_l} \int^\infty \{\hat H_i^{(-l)}(t; \theta) - H_{i}(t; \theta)\}dQ_{i}(t)\Bigg| \cO_{\cI_{-l}}\right]\right)\\
    &~~~+\frac{1}{n\delta^2}\Var\left(\E\left[\sum_{i \in \cI_l} \int^\infty \{\hat H_i^{(-l)}(t; \theta) - H_{i}(t; \theta)\}dQ_{i}(t)\Bigg| \cO_{\cI_{-l}}\right]\right)\\
    &=\frac{1}{n\delta^2}\E\left\{\E\left(\left[\sum_{i \in \cI_l} \int^\infty \{\hat H_i^{(-l)}(t; \theta) - H_{i}(t; \theta)\}dQ_{i}(t)\right]^2\Bigg| \cO_{\cI_{-l}}\right)\right\}\\
    &=\frac{1}{n\delta^2}\E\left\{\sum_{i \in \cI_l} \E\left(\left[\int^\infty \{\hat H_i^{(-l)}(t; \theta) - H_{i}(t; \theta)\}dQ_{i}(t)\right]^2\Bigg| \cO_{\cI_{-l}}\right)\right\}\\
    &=\frac{1}{\delta^2L}\E\left\{\E\left(\left[\int^\infty \{\hat H(t; \theta) - H(t; \theta)\}dQ(t)\right]^2\Bigg| \cO_{\cI_{-l}}\right)\right\}\\
    &=\frac{1}{\delta^2L}\E_{\cO_{\cI_{-l}}}\left\{\E_{\cO_{\cI_{l}}}\left(\left[\int^\infty \{\hat H(t; \theta) - H(t; \theta)\}dQ(t)\right]^2\right)\right\}\\
    &\leq \frac{1}{\delta^2L}\E_{\cO_{\cI_{-l}}}\left(\E_{\cO_{\cI_{l}}}\left[\sup_t |\hat H(t; \theta) - H(t; \theta)|^2 V^\infty(Q)^2\right]\right)\\
    &\leq \frac{1}{\delta^2L}\E\left\{\|\hat H(\theta) - H(\theta)\|_{\sup, 2}^2\right\}\|V^\infty(Q)^2\|_\infty\\
    &\leq \frac{1}{\delta^2L}\E\left\{\|\hat H(\theta) - H(\theta)\|_{\sup, 2}^2\right\}\|V^\infty(Q)\|_\infty^2\\
    &= o(1),
\end{align*}
}
by Assumption \ref{assump:uniformcons2}. $\cT_3$ follows the same argument as that of $\cT_2$ immediately. The arguments on $\cT_4$, $\cT_5$, $\cT_6$ are the same as in the proof of Theorem \ref{thm:mdr} because their definitions do not change with whether cross fitting is employed. Therefore consistency and asymptotic normality follow immediately.

\bibliographystyle{imsart-nameyear}
\bibliography{ref,lpci_ref,left_truncation}

\end{document}

%% file: notation.tex
\newtheorem{thm}{Theorem}

\newtheorem{assump}{Assumption}

\theoremstyle{remark}

\newtheorem{exm}{Example}

\def\beq{\begin{equation}} 
\def\eeq{\end{equation}}
\def\beqn{\begin{eqnarray*}}
\def\eeqn{\end{eqnarray*}}
\def\Bitem{\begin{itemize}\setlength{\itemsep}{.2in}}
\def\bitem{\begin{itemize}\setlength{\itemsep}{.05in}}
\def\eitem{\end{itemize}}
\def\Benum{\begin{enumerate}\setlength{\itemsep}{.2in}}
\def\benum{\begin{enumerate}\setlength{\itemsep}{.05in}}
\def\eenum{\end{enumerate}}
\def\bmult{\begin{multline*}}
\def\emult{\end{multline*}}
\def\bcenter{\begin{center}}
\def\ecenter{\end{center}}
\def\bframe{\begin{frame}}
\def\eframe{\end{frame}}

\def\cI{\mathcal{I}}

\def\cN{\mathcal{N}}
\def\cO{\mathcal{O}}
\def\cP{\mathcal{P}}

\def\cT{\mathcal{T}}

\newcommand{\E}{\operatorname{\mathbb{E}}}
\renewcommand{\P}{\operatorname{\mathbb{P}}}
\newcommand{\Var}{\operatorname{Var}}

\def\1{\mathbbm{1}}